\documentclass[11pt, a4paper]{article}
\usepackage{graphicx}
\usepackage{mathrsfs}
\usepackage{CJK}
\usepackage{bbm}
\usepackage{stmaryrd}
\usepackage{amsmath}
\usepackage{cases}
\usepackage{amsfonts}
\usepackage{latexsym,amssymb,amsmath,mathrsfs,amsbsy, amsthm}
\usepackage[usenames]{color}
\usepackage{xspace,colortbl}
\usepackage[dvipdfm,
pdfstartview=FitH, CJKbookmarks=true, bookmarksnumbered=true,
bookmarksopen=true,
citecolor=blue 
,colorlinks=true
,pdfborder=1 
,pdfstartpage=2 
,pdfstartview=Fit]{hyperref} 
\allowdisplaybreaks

\newtheorem{thm}{Theorem}[section]
\newtheorem{lem}[thm]{Lemma}

\newtheorem{rem}{Remark}

\newcommand{\nn}{\nonumber}

\numberwithin{equation}{section}
\begin{document}
\title{\bf\Large On Vortex Solutions of the Landau-Lifshitz Equations}

\author{Ruiqi JIANG\thanks{Supported by NSFC, Grant No. 10990013}}

\date{}
\maketitle

\begin{abstract}
We study the Landau-Lifshitz equation of ferromagnetism on $\mathbb{R}^2$ with an easy-axis anisotropy. We give the necessary condition for the existence of the finite energy vortex solutions
and show the behaviors of the solutions.
\end{abstract}


\section{Introduction}

The Landau-Lifshitz equation describes the magnetization phenomenon in ferromagnetic medium. Here we
consider easy-axis anisotropic case in two space dimensions
\begin{eqnarray}\label{eqn:1.1}
  \frac{\partial u}{\partial t}=u \times (\Delta u +\lambda u_3 e_3)
\end{eqnarray}
where $u=(u_1,u_2,u_3):\mathbb{R}^2 \to \mathbb{S}^2 \subset \mathbb{R}^3$, $e_3=(0,0,1)$ denotes the north pole.

Although the local and global existence of smooth solution with small data has been established
(see\cite{DW},\cite{BIT1},\cite{CSU} and the references therein), the global existence with large initial data remains open. Recently several authors(\cite{BIT2}) have obtained the global equivariant solutions with energy less than $4\pi$. While Lin and Wei (\cite{LW}) constructed traveling wave
solutions to equation(\ref{eqn:1.1}) with $\lambda<0$, we would like to study some topologically nontrivial, periodic solutions known as vortex or vortex-like solutions.

In this paper, we seek for solutions which are equivariant with respect to the $S^1=O(2)$ actions on
both $\mathbb{R}^2$ and $\mathbb{S}^2$. Specifically, we look for a solution of the following form
$$u(r,\theta)=(\sin h(r) \sin(m\theta+\omega t+\theta_0) ,\,\sin h(r)\cos(m\theta+\omega t+\theta_0),\,\cos h(r))$$
where $(r,\theta)$ denotes the polar coordinates in $\mathbb{R}^2$, $m\in \mathbb{Z}$ is a
topological degree (known as vortex degree in physics), $\omega \in \mathbb{R}$ is the angular velocity and $\theta_0$ is the initial phase. By direction calculations, equations (\ref{eqn:1.1}) reduces to an ordinary differential equation (ODE) of $h(r)$
\begin{eqnarray}\label{eqn:1.2}
h''(r)+\frac{1}{r}h'(r)-\frac{m^2}{r^2}\sin h(r) \cos h(r)=\lambda \sin h(r) \cos h(r)+ \omega \sin h(r)
\end{eqnarray}
These time periodic solutions are called magnetic vortices or vortex solutions, which play a important role in the geometry and topology of the flow.

Then, for $u=(u_1,u_2,u_3)\in \mathbb{R}^3$, the energy is simply $$E(u)=\frac{1}{2}\int_{\mathbb{R}^2}{|\nabla u|^2dx}$$ where
$$|\nabla u|^2=\sum_{i=1}^{3}|\nabla u_i|^2$$
For this particular form of solutions, the energy $E$ reduces to a functional $J$ on the function $h$ as follows.(We omit the factor $\pi$ in the integrals)
$$J(h)=\int_0^{+\infty}{(h'^2+\frac{m^2}{r^2}\sin^2 h)rdr}$$

In order to find smooth solutions to equation (\ref{eqn:1.1}), let's consider the following initial
value problem to ODE (\ref{eqn:1.2})
\begin{eqnarray}\label{eqn:1.3}
h(0)=0,\quad h^{(|m|)}(0)=a,
\end{eqnarray}
where $a\in \mathbb{R}$ and $h^{(|m|)}$ denotes the $|m|$-order derivative of $h(r)$.

For $m \in \mathbb{Z}\setminus \{0\}$, we have prove existence, uniqueness and continuous dependence on the initial data to the problem (\ref{eqn:1.2})-(\ref{eqn:1.3}) in \cite{J}. Since the vortex solutions tend to a fixed point on the sphere $\mathbb{S}^2$ as $|x|\to \infty$, it reduces to the following boundary condition of $h(r)$,
\begin{eqnarray}\label{eqn:1.4}
\lim_{r \to \infty} h(r)=k\pi,\quad \text{for some}\ k \in \mathbb{Z}.
\end{eqnarray}
In this paper, we say a solution to the problem (\ref{eqn:1.2})-(\ref{eqn:1.3}) is a vortex solution if it satisfies the boundary condition (\ref{eqn:1.4}).

For $\lambda =0, \omega=0$, there exist explicit solutions to the problem (\ref{eqn:1.2})-(\ref{eqn:1.3}), known as Belavin-Polyakov instantons(\cite{K1}). These solutions, which correspond to a family of harmonic maps from $\mathbb{R}^2$ to $\mathbb{S}^2$, have the following expression:
\begin{eqnarray}
h(r)=2\arctan[\frac{a}{2|m|!}r^{|m|}], \quad |m|>0.
\end{eqnarray}

For $\lambda =0$, $\omega \neq 0$, Kollar (\cite{K}) established the nonexistence of finite-energy nontrivial votex solutions to the problem (\ref{eqn:1.2})-(\ref{eqn:1.3}) with $m \neq 0$.

For $\lambda >0$, $\omega>0$, the existence of finite-energy nontrivial vortex solutions to (\ref{eqn:1.2})-(\ref{eqn:1.3}) has been studied, see \cite{G_S} and \cite{J}.

In this paper, we get the necessary condition about the existence of the finite-energy nontrivial vortex solutions by ruling out all the other possibilities and show the behaviors of the solutions. Now
let's state our main result.
\begin{thm}\label{thm:1}
For $\lambda,\, \omega \in \mathbb{R},\, m \in \mathbb{Z}\setminus \{0\}$, if there exists a finite-energy vortex solution $h(r)$ to the problem
(\ref{eqn:1.2})-(\ref{eqn:1.3}) with $a \neq 0$, then exactly one of the following holds:\\
(i).~~$\lambda=\omega=0$, $h(r)=2\arctan[\frac{a}{2|m|!}r^{|m|}]$.\\
(ii).~~$\lambda=\omega \neq 0$  and $h(r)$ converges to $(2l+1)\pi$ as $ r\to \infty$ for some $l \in \mathbb{Z}$.\\
(iii).~$\lambda=-\omega \neq 0$ and $h(r)$ converges to $2l\pi$ as $ r\to \infty$ for some $l \in \mathbb{Z}$.\\
(iv).~~$0<\omega < \lambda$ and $h(r)$ converges to $(2l+1)\pi$ exponentially as $r \to \infty$ for some $l \in \mathbb{Z}$.\\
(v).~~$-\lambda < \omega <0$ and $h(r)$ converges to $2l\pi$ exponentially as $r \to \infty$ for some $l \in \mathbb{Z}$.
\end{thm}.

\begin{rem}\label{rem:1}
Since it's not easy to characterize the behavior of solutions to (\ref{eqn:1.2}) in case (ii) and (iii), especially its rate of convergence as $r \to \infty$, we can't rule out the case (ii) and (iii) by Pohozaev identity. We believe that case (ii) and (iii) will never happen.
\end{rem}

\noindent \textbf{Convention:} For convenience, we always assume that $m>0$ without further comment.

\section{Nonexistence of the vortices solutions}
In this section we will establish two theorems about nonexistence of the vortices solutions.

Let's consider following equation
\begin{eqnarray}\label{eqn:2.1}
h''(r)+\frac{1}{r}h'(r)-\frac{m^2}{r^2}\sin h(r) \cos h(r)=g(h(r)),\quad r \in (0,+\infty).
\end{eqnarray}
and corresponding initial values
\begin{eqnarray}\label{eqn:2.2}
h(0)=0,\quad h^{(m)}(0)=a.
\end{eqnarray}
where $g(\cdot) \in C^{\infty}(\mathbb{R})$ and $a \in \mathbb{R}$.

\begin{lem}\label{lem:2.1}
Suppose nonconstant function $h(r)$ satisfies equation (\ref{eqn:2.1}) and $g(k_0\pi)=0$, $g'(k_0\pi)>0$ for some $k_0 \in \mathbb{Z}$. If $\lim \limits_{r \to \infty}h(r)=k_0\pi$, there exists $R_0>0$ such that exactly one of the following holds:\\
(1)$h(r)$ decreases to $k_0\pi$ exponentially as $r \to \infty$ on the interval $[R_0,\, \infty)$.\\
(2)$h(r)$ increases to $k_0\pi$ exponentially as $r \to \infty$ on the interval $[R_0,\, \infty)$.
\end{lem}

\begin{proof}
Let $\widetilde{h}(r)=h(r)-k_0\pi$, then $\widetilde{h}(r)$ satisfies the following equation
\begin{eqnarray}\label{eqn:2.3}
\widetilde{h}''(r)+\frac{1}{r}\widetilde{h}'(r)-\frac{m^2}{r^2}\sin \widetilde{h}(r) \cos \widetilde{h}(r)=\widetilde{g}(\widetilde{h}(r)),\quad r \in (0,+\infty).
\end{eqnarray}
where $$\widetilde{g}(\widetilde{h}(r))=g(\widetilde{h}(r)+k_0\pi).$$

Now there hold that $$\lim \limits_{r \to \infty}\widetilde{h}(r)=0, \quad \widetilde{g}(0)=0\quad \text{and} \quad  \widetilde{g}'(0)>0$$

For $\widetilde{g}'(0)>0$, there exist $\delta_0>0$ and $r_0>0$ such that for any $(r,x) \in
[r_0,\infty)\times[-\delta_0,\delta_0]$ there hold
\begin{eqnarray}\label{eqn:2.4}
\frac{m^2}{r^2}\sin x \cos x + \widetilde{g}(x)>\frac{1}{2}\widetilde{g}'(0)x, \quad x>0
\end{eqnarray}
and
\begin{eqnarray}\label{eqn:2.5}
\frac{m^2}{r^2}\sin x \cos x + \widetilde{g}(x)\leq \frac{1}{2}\widetilde{g}'(0)x, \quad x\leq0
\end{eqnarray}

For $\lim \limits_{r \to \infty}\widetilde{h}(r)=0$ and $\widetilde{h}(r)$ is not a constant function, there exists $R_0>r_0$ such that $$\widetilde{h}(R_0)\neq 0,\quad \text{and} \quad |\widetilde{h}(r)|<\delta_0,\quad r>R_0.$$

{\large \textbf{Case I.}}: If $\widetilde{h}(R_0)>0$, we claim that
$$\widetilde{h}(r)>0\quad  \text{and}\quad \widetilde{h}'(r)<0,\quad r \in [R_0, \infty).$$

If there exists $r^* \in [R_0,\infty)$ such that $\widetilde{h}(r^*)\leq 0$, combining $\widetilde{h}(R_0)>0$ and $\lim \limits_{r \to \infty}\widetilde{h}(r)=0$, it's easy to see that $\widetilde{h}(r)$ has a minimum point
$\widehat{r}\geq r^*$ such that
$$\widetilde{h}(\widehat{r})\leq 0,\quad \widetilde{h}'(\widehat{r})=0\quad \text{and}\quad \widetilde{h}''(\widehat{r})\geq 0.$$

Since $\widetilde{h}(\widehat{r})=0,\ \widetilde{h}'(\widehat{r})=0$ leads to $\widetilde{h}(r)\equiv 0$ which contradicts the fact that $h(r)=\widetilde{h}(r)+k_0\pi$ is not a constant funciton , we have $\widetilde{h}(\widehat{r})<0$. However, by the equation (\ref{eqn:2.3}) and inequality (\ref{eqn:2.5}), we have
$$0\leq \widetilde{h}''(\widehat{r})=\frac{m^2}{2\widehat{r}^2}\sin 2\widetilde{h}(\widehat{r})+\widetilde{g}(\widetilde{h}(\widehat{r}))\leq \frac{1}{2}\widetilde{g}'(0)\widetilde{h}(\widehat{r})<0,$$
there exists a contradiction. So we obtain that $$\widetilde{h}(r)>0,\quad r\in [R_0,\infty).$$

Next we will prove $$\widetilde{h}'(r)<0,\quad r\in [R_0,\infty).$$

From equation (\ref{eqn:2.3}), we have that, for $R_0\leq s \leq r<+\infty$
\begin{eqnarray}\label{eqn:2.6}
r\widetilde{h}'(r)=s\widetilde{h}'(s)+\int_s^r{[\frac{m^2}{2t^2}\sin 2\widetilde{h}(t)+\widetilde{g}(\widetilde{h}(t))]tdt}
\end{eqnarray}
and
\begin{eqnarray}\label{eqn:2.7}
\widetilde{h}(r)=\widetilde{h}(s)+\int_s^r{\widetilde{h}'(t)dt}.
\end{eqnarray}

Combining $0<\widetilde{h}(r)<\delta_0,\ r\in [R_0,\infty)$ and (\ref{eqn:2.4}), we get that $r\widetilde{h}'(r)$ increases monotonically on the interval $[R_0,\infty)$. Since $\widetilde{h}(r)$ is bounded on the interval $[R_0,\infty)$, (\ref{eqn:2.7}) implies that there exist $C>0$ and $r_k\to \infty$ such that $$\widetilde{h}'(r_k)\leq \frac{C}{r_k}.$$
Thus, there holds that $\lim \limits_{r \to \infty} r\widetilde{h}'(r)<\infty$. Moreover $\lim \limits_{r \to \infty} r\widetilde{h}'(r)=0$, otherwise $\widetilde{h}$ will be unbounded on the interval $[R_0,\infty)$. Let $r\to \infty$ and replace $s$ by $r$, (\ref{eqn:2.6}) reduces to
$$r\widetilde{h}'(r)=-\int_r^{\infty}{[\frac{m^2}{2t^2}\sin 2\widetilde{h}(t)+\widetilde{g}(\widetilde{h}(t))]tdt}$$

By (\ref{eqn:2.4}) and $0<\widetilde{h}(r)<\delta_0,\ r\in [R_0,\infty)$, we deduce that
$$\widetilde{h}'(r)<0,\quad r\in [R_0,\infty).$$

Now we are in position to show $\widetilde{h}(r)$ decreases to 0 exponentially as $r \to \infty$.

Let $f(r)=be^{-\epsilon r}$, where $b,\,\epsilon>0$ will be determined later and denote $$\beta(r)=f(r)-\widetilde{h}(r),$$
it follows that
$$\beta''+\frac{1}{r}\beta'=(\epsilon^2-\frac{\epsilon}{r})f-\frac{m^2}{2r^2}\sin 2\widetilde{h}-\widetilde{g}(\widetilde{h})$$

By (\ref{eqn:2.4}) and $0<\widetilde{h}(r)<\delta_0,\ r\in [R_0,\infty)$, we obtain that
$$\beta''(r)+\frac{1}{r}\beta'(r) < \epsilon^2f(r)-\frac{1}{2}\widetilde{g}'(0)\widetilde{h}(r),\quad r \in [R_0,\infty)$$
Choose $\epsilon=\sqrt {\frac{1}{2}\widetilde{g}'(0)}$, it deduces that
$$\beta''(r)+\frac{1}{r}\beta'(r) -\frac{1}{2}\widetilde{g}'(0)\beta(r)<0,\quad r \in [R_0,\infty)$$

After choosing $b=b_0$ such that $$\beta(R_0)=b_0e^{-\sqrt {\frac{1}{2}\widetilde{g}'(0)} R_0}-\widetilde{h}(R_0)>0,$$
by maximum principle we have $\beta(r)\geq 0,\  r \in [R_0,\infty)$, i.e.
$$0<\widetilde{h}(r)=h(r)-k_0\pi\leq b_0e^{-\sqrt {\frac{1}{2}\widetilde{g}'(0)} r},\quad r \in [R_0,\infty)$$

By above argument, we know $h(r)$ decreases to $k_0\pi$ exponentially as $r \to \infty$ on the interval $[R_0, \infty)$.

{\large \textbf{Case II.}}: If $\widetilde{h}(R_0)<0$, following the similar argument, we have that  $h(r)$ increases to $k_0\pi$ exponentially as $r \to \infty$ on the interval $[R_0, \infty)$.

\end{proof}

\begin{thm}\label{thm:2.2}
If $g(k_0\pi)=0$, $g'(k_0\pi)>0$ for some $k_0 \in \mathbb{Z}$ and $G(x)=-\int_x^{k_0\pi}{g(t)dt}$ doesn't change sign on $\mathbb{R}$, there is no solution to the problem (\ref{eqn:2.1})-(\ref{eqn:2.2}) with $a \neq 0$ which satisfies
$$\lim \limits_{r \to \infty}h(r)=k_0\pi.$$
\end{thm}

\begin{proof}
Multiplying both sides of equation (\ref{eqn:2.1}) by $r^2h'(r)$ and integrating on the interval
$[0,r]$ yields the Pohozaev identity
\begin{eqnarray}
(rh'(r))^2&=&m^2\sin^2 h(r)+2\int_0^r{g(h(t))h'(t)t^2dt}\nn\\
&=&m^2\sin^2 h(r)+2G(h(r))r^2-4\int_0^r{G(h(t))tdt}\nn
\end{eqnarray}

If there is a solution $h(r)$ to the problem (\ref{eqn:2.1})-(\ref{eqn:2.2}) which satisfies $\lim \limits_{r \to \infty}h(r)=k_0\pi$,
by Lemma(\ref{lem:2.1}), we have that
$$\lim_{r \to \infty}rh'(r)=0\quad \text{and}\quad \lim_{r \to \infty}r^2(h(r)-k_0\pi)^2=0.$$

Since $G(x)=g'(k_0\pi)(x-k_0\pi)^2 + o((x-k_0\pi)^2)$, let $r\to \infty$, then Pohozaev identity reduces to $$\int_0^\infty{G(h(t))tdt}=0.$$
However, the fact that $G(x)$ doesn't change sign on $\mathbb{R}$ implies that $h(r) \equiv k_0\pi$.
For $h(0)=0$, we obtain $k_0=0$ and $a=0$. The theorem is proved.

\end{proof}

 The following theorem is an analogue to the results of Kollar (\cite{K}) for infinite energy of the oscillation solutions. For the reader's convenience, we give its proof.

\begin{thm}\label{thm:2.3}
Suppose nonconstant function $h(r)$ satisfies equation (\ref{eqn:2.1}) and $g(k_0\pi)=0$, $g'(k_0\pi)<0$ for some $k_0 \in \mathbb{Z}$. If $\lim \limits_{r \to \infty}h(r)=k_0\pi$,
then there exists $R_0>0$ such that $h(r)$ oscillates around $k_0\pi$ on the interval $(R_0,\infty)$
and $$\int_{R_0}^\infty{[h'^2+\frac{m^2}{r^2}\sin^2 h]rdr}=+\infty.$$
\end{thm}
\begin{proof}
We will complete the proof through four steps.

{\large \textbf{Step one}}: Let $\widetilde{h}(r)=h(r)-k_0\pi$, then $\widetilde{h}(r)$ satisfies the following equation
\begin{eqnarray}\label{eqn:2.8}
\widetilde{h}''(r)+\frac{1}{r}\widetilde{h}'(r)-\frac{m^2}{r^2}\sin \widetilde{h}(r) \cos \widetilde{h}(r)=\widetilde{g}(\widetilde{h}(r)),\quad r \in (0,+\infty).
\end{eqnarray}
where $\widetilde{g}(\widetilde{h}(r))=g(\widetilde{h}(r)+k_0\pi)$.

Now there hold that $$\lim \limits_{r \to \infty}\widetilde{h}(r)=0, \quad \widetilde{g}(0)=0\quad \text{and} \quad  \widetilde{g}'(0)<0$$

For $\widetilde{g}'(0)<0$, there exist $\delta_0>0$ and $r_0>0$ such that for any $(r,x) \in
[r_0,\infty)\times[-\delta_0,\delta_0]$ there hold
\begin{eqnarray}\label{eqn:2.9}
\frac{3}{2}\widetilde{g}'(0) \leq \frac{m^2}{r^2}\frac{\sin 2x}{2x} + \frac{\widetilde{g}(x)}{x}\leq \frac{1}{2}\widetilde{g}'(0).
\end{eqnarray}

For $\lim \limits_{r \to \infty}\widetilde{h}(r)=0$, there exists $R_0>r_0$ such that
\begin{eqnarray}\label{eqn:2.10}
|\widetilde{h}(r)|<\delta_0,\quad r\in [R_0,\,\infty).
\end{eqnarray}

\medskip
{\large \textbf{Step two}}: We would like to prove that $\widetilde{h}(r)$ oscillate around $0$ on the interval $(R_0,\infty)$ which is equivalent to $h(r)$ oscillating around $k_0\pi$ on the same interval.

Let $H(r)=r^{\frac{1}{2}}\widetilde{h}(r)$, then $H(r)$ satisfies the following equation
\begin{eqnarray}\label{eqn:2.11}
H''(r)+C(r)H(r)=0, \quad r \in (R_0,\infty)
\end{eqnarray}
where $$C(r)=\frac{1}{4r^2}-(\frac{m^2}{r^2}\frac{\sin 2\widetilde{h}}{2\widetilde{h}} + \frac{\widetilde{g}(\widetilde{h})}{\widetilde{h}})$$

By (\ref{eqn:2.9}) and (\ref{eqn:2.10}), we have that
$$C(r)>-(\frac{m^2}{r^2}\frac{\sin 2\widetilde{h}}{2\widetilde{h}} + \frac{\widetilde{g}(\widetilde{h})}{\widetilde{h}})\geq -\frac{1}{2}\widetilde{g}'(0)>0$$

Comparing equation (\ref{eqn:2.11}) with the following equation
$$u''(r)+l^2u(r)=0$$
where $l=\sqrt{-\frac{1}{2}\widetilde{g}'(0)}$, whose solutions have period $\frac{2\pi}{l}$ and
infinite zero points in $\mathbb{R}$.

By Sturm-Liouville theorem, we get $H(r)$ has at least one zero point between any two adjacent zero
points of $u(r)$ on the interval $(R_0,\infty)$. That's to say, nonconstant function $\widetilde{h}(r)=r^{-\frac{1}{2}}H(r)$ must oscillates around zero.

\medskip
{\large \textbf{Step three}}: In this step we will show that $\widetilde{h}(r)$ has similar monotonicity as function $\sin r$ or $\cos r$.

Let $R_0\leq a_1 <a_2<\infty$ be two adjacent zero points of $\widetilde{h}(r)$, i.e $\widetilde{h}(a_1)=\widetilde{h}(a_2)=0$. For $\widetilde{h}(r)$ is not a constant function, we deduce that $\widetilde{h}'(a_1)\neq 0$. Without loss of generality, we assume $\widetilde{h}'(a_1)>0$. So there holds that
$$0<\widetilde{h}(r)< \delta_0,\quad r \in (a_1,a_2).$$

By (\ref{eqn:2.8}) and (\ref{eqn:2.9}), we obtain
$$\widetilde{h}''(r)+\frac{1}{r}\widetilde{h}'(r)=\frac{m^2}{2r^2}\sin 2\widetilde{h}(r) +\widetilde{g}(\widetilde{h}(r))\leq \frac{1}{2}\widetilde{g}'(0)\widetilde{h}(r)<0,\quad r \in (a_1,a_2)$$
Above inequality says that $\widetilde{h}(r)$ doesn't have local minima in the interval $(a_1,a_2)$. It implies
that $\widetilde{h}(r)$ must increase monotonically from 0 to unique local maxima, then decrease monotonically to zero on the interval $[a_1,a_2]$. For the case $\widetilde{h}'(a_1)<0$, following the similar argument, there holds that $\widetilde{h}(r)$ must decrease monotonically from 0 to unique local minima, then increase monotonically to zero on the interval $[a_1,a_2]$.

Thus, it's convenient to introduce $\{a_i\}$ $\{M_i\}$ $\{L_i\}$, the increasing infinite sequences of zero points of $\widetilde{h}(r)$, the sequences of local maxima and local minima of $\widetilde{h}(r)$, respectively. By neglecting first few terms, we may assume $$R_0<a_1<M_1<a_2<L_1<a_3<\cdot\cdot\cdot\cdot.$$

\medskip
{\large \textbf{Step four}}:

(i) Let $e(r)=r(\widetilde{h}'(r))^2$, $r\in [a_{2l-1},M_l]$ for $l \in \mathbb{N}$, from the conclusion of step three, we have that $$\widetilde{h}'(r)>0\ \text{and}\ \widetilde{h}(r)>0,\quad r\in(a_{2l-1},M_l).$$

By direct calculation, for $r\in(a_{2l-1},M_l)$, we have that
\begin{eqnarray}
\frac{de(r)}{dr}
&=&(\widetilde{h}'(r))^2+2r\widetilde{h}'(r)\widetilde{h}''(r)\nn\\
&=&(\widetilde{h}'(r))^2+2r\widetilde{h}'(r)[-\frac{1}{r}\widetilde{h}'(r)+\frac{m^2}{2r^2}\sin 2\widetilde{h}+\widetilde{g}(\widetilde{h})]\nn\\
&=&-(\widetilde{h}'(r))^2+2r\widetilde{h}'(r)[\frac{m^2}{2r^2}\sin 2\widetilde{h}+\widetilde{g}(\widetilde{h})]\nn\\
&>&-\frac{e(r)}{r}+3r\widetilde{h}'(r)\widetilde{h}(r)\widetilde{g}'(0)>-\frac{e(r)}{r}
\end{eqnarray}

Integrating above inequality on $[a_{2l-1},r)$, it follows that
$$e(r)\geq e(a_{2l-1})\frac{a_{2l-1}}{r},\quad r \in [a_{2l-1},M_l].$$
Integrating above inequality over $r \in [a_{2l-1},M_l]$ yields that
$$\int_{a_{2l-1}}^{M_l}{e(r)dr}\geq e(a_{2l-1})a_{2l-1}\log {\frac{M_l}{a_{2l-1}}},$$
i.e.
\begin{eqnarray}\label{eqn:2.14}
\int_{a_{2l-1}}^{M_l}{(\widetilde{h}'(r))^2 rdr}\geq (a_{2l-1}\widetilde{h}'(a_{2l-1}))^2\log {\frac{M_l}{a_{2l-1}}}.
\end{eqnarray}

Following the similar argument, we obtain,
\begin{eqnarray}\label{eqn:2.15}
\int_{a_{2l}}^{L_l}{(\widetilde{h}'(r))^2 rdr}\geq (a_{2l}\widetilde{h}'(a_{2l}))^2\log {\frac{L_l}{a_{2l}}}.
\end{eqnarray}

(ii) By (\ref{eqn:2.8}) and (\ref{eqn:2.9}), we have that, for $k \in \mathbb{N}$,
\begin{eqnarray}
 (r\widetilde{h}'(r))^2|_{a_k}^{a_{k+1}}
&=&m^2\sin^2 \widetilde{h}(r)|_{a_k}^{a_{k+1}}+2\int_{a_k}^{a_{k+1}}{\widetilde{g}(\widetilde{h})\widetilde{h}'r^2dr}\nn\\
&=&2\int_{a_k}^{a_{k+1}}{\frac{\widetilde{g}(\widetilde{h})}{\widetilde{h}}\widetilde{h}\widetilde{h}'r^2dr}\nn\\
&\sim& 2\widetilde{g}'(0)\int_{a_k}^{a_{k+1}}{\widetilde{h}\widetilde{h}'r^2dr}=-2\widetilde{g}'(0)\int_{a_k}^{a_{k+1}}{\widetilde{h}^2 rdr}>0
\end{eqnarray}
Immediately it yields that, for $k \in \mathbb{N}$,
\begin{eqnarray}\label{eqn:2.17}
(a_{k+1}\widetilde{h}'(a_{k+1}))^2>(a_k\widetilde{h}'(a_k))^2.
\end{eqnarray}

(iii) Let $\theta(r)=\arctan \frac{\widetilde{h}'(r)}{\alpha \widetilde{h}(r)}$, $r \in (a_k,a_{k+1})$, where $\alpha=\sqrt{-\widetilde{g}'(0)}$ and $k \in \mathbb{N}$.

By direct calculation,
$$\theta'(r)=-\alpha+[-\frac{1}{r}\frac{\alpha \widetilde{h}\widetilde{h}'}{(\alpha \widetilde{h})^2+\widetilde{h}'^2}+\frac{m^2}{2r^2}
\frac{\alpha \widetilde{h}\sin 2\widetilde{h}}{(\alpha \widetilde{h})^2+\widetilde{h}'^2} + \frac{\alpha \widetilde{h}(\widetilde{g}(\widetilde{h})+\alpha^2\widetilde{h})}{(\alpha \widetilde{h})^2+\widetilde{h}'^2}],$$
it yields that
$$|\theta'(r)+\alpha|\leq \frac{1}{r}+\frac{2m^2}{r^2}+\alpha|\widetilde{h}|\frac{O(\widetilde{h}^2)}{(\alpha \widetilde{h})^2+\widetilde{h}'^2}$$

So $\forall \epsilon>0$, there exists $K \in \mathbb{N}$ such that for $k \geq K$, there holds
$$|\theta'(r)+\alpha|\leq \epsilon$$

For $l\in \mathbb{N}$, there holds that
$$\lim_{r \to a_{2l-1}^+} \theta(r)=\frac{\pi}{2},\quad \lim_{r \to a_{2l}^-} \theta(r)=-\frac{\pi}{2},\quad \theta(M_l)=0;$$
$$\lim_{r \to a_{2l}^+} \theta(r)=\frac{\pi}{2},\quad \lim_{r \to a_{2l+1}^-} \theta(r)=-\frac{\pi}{2},\quad \theta(L_l)=0.$$

Integrating $\theta'(r)$ on $[a_{2l-1},M_l]$ and $[M_l,a_{2l}]$ yield that
$$(-\alpha -\epsilon)(M_l-a_{2l-1})\leq -\frac{\pi}{2}=\int_{a_{2l-1}}^{M_l}{\theta'(r)dr} \leq
(-\alpha +\epsilon)(M_l-a_{2l-1})$$
$$(-\alpha -\epsilon)(a_{2l}-M_l)\leq -\frac{\pi}{2}=\int_{M_l}^{a_{2l}}{\theta'(r)dr} \leq
(-\alpha +\epsilon)(a_{2l}-M_l)$$

Combining above two inequalities, we have that
$$|M_l-\frac{a_{2l-1}+a_{2l}}{2}|\leq \pi \epsilon\quad \text{and}\quad |\frac{a_{2l}-a_{2l-1}}{\pi}-\alpha|\leq \epsilon$$
It's easy to verify that there exist $\mu \in (0,1)$ such that
\begin{eqnarray}\label{eqn:2.18}
\log \frac{M_l}{a_{2l-1}}\geq \mu \log \frac{a_{2l}}{a_{2l-1}},\quad l >K.
\end{eqnarray}
By the similar argument, we also obtain
\begin{eqnarray}\label{eqn:2.19}
\log \frac{L_l}{a_{2l}}\geq \mu \log \frac{a_{2l+1}}{a_{2l}},\quad l >K.
\end{eqnarray}

(iv) Combining (\ref{eqn:2.14}), (\ref{eqn:2.15}), (\ref{eqn:2.17}), (\ref{eqn:2.18}) and (\ref{eqn:2.19}), we have
\begin{eqnarray}
\int_{a_1}^\infty{(\widetilde{h}'(r))^2rdr}
&=& \sum_{k=1}^\infty \int_{a_k}^{a_{k+1}}{(\widetilde{h}'(r))^2rdr}\nn\\
&\geq& \sum_{l=K}^\infty [\int_{a_{2l-1}}^{M_l}{(\widetilde{h}'(r))^2rdr}+\int_{a_{2l}}^{L_l}{(\widetilde{h}'(r))^2rdr}\nn\\
&\geq& ((a_{2K-1}\widetilde{h}'(a_{2K-1}))^2 \sum_{l=K}^\infty [\log \frac{M_l}{a_{2l-1}}+\log \frac{L_l}{a_{2l}}]\nn\\
&\geq& \mu ((a_{2K-1}\widetilde{h}'(a_{2K-1}))^2 \sum_{l=K}^\infty [\log \frac{a_{2l}}{a_{2l-1}}+\log \frac{a_{2l+1}}{a_{2l}}]=\infty
\end{eqnarray}

Immediately it follows that
$$\int_{R_0}^\infty{[h'^2+\frac{m^2}{r^2}\sin^2 h]rdr}=\int_{R_0}^\infty{[\widetilde{h}'^2+\frac{m^2}{r^2}\sin^2 \widetilde{h}]rdr} \geq \int_{a_1}^\infty{(\widetilde{h}'(r))^2rdr}=\infty$$

\end{proof}

\section{The Proof of Theorem \ref{thm:1}}
In this section, we will employ theorem (\ref{thm:2.2}) and (\ref{thm:2.3}) to prove our main theorem by ruling out the cases that vortex solutions with finite energy don't exist.

\medskip
\noindent{\bf\textit{Proof of Theorem \ref{thm:1}}}.

Let $g(h)=\lambda \sin h \cos h+\omega \sin h$, then $$g'(h)=\lambda \cos 2h + \omega \cos h.$$
Define $$G(x,k)=-\int_x^{k\pi}{g(t)dt}=\frac{\lambda}{2}\sin^2 x +\omega[(-1)^{|k|}-\cos x].$$

\medskip
\noindent {\large \bf {Case one}}: $\lambda=\omega=0$.

we can get solutions to the problem (\ref{eqn:1.2})-(\ref{eqn:1.3}) which has following explicit expression: $$h(r)=2\arctan[\frac{a}{2m!}r^{m}],\quad m>0.$$

\noindent {\large \bf {Case two}}: $\lambda=0,\,\omega \neq0$.

Now we have $g'(k\pi)=\omega (-1)^{|k|}\neq 0$ for $k \in \mathbb{Z}$.

(i) When $g'(k\pi)<0$, by theorem (\ref{thm:2.3}), there don't exist vortex solutions with finite energy which satisfy $\lim_{r \to \infty} h(r)=k\pi$.

(ii) When $g'(k\pi)>0$, for $G(x,k)=\omega[(-1)^{|k|}-\cos x]$ doesn't change sign on $\mathbb{R}$, by theorem (\ref{thm:2.2}), there don't exist vortex solutions which satisfy $\lim_{r \to \infty} h(r)=k\pi$.

From above argument, we know there don't admit any vortex solutions with finite energy in this case.

\medskip
\noindent {\large \bf {Case three}}: $\lambda \neq 0,\,\omega =0$.

Now we have $g'(k\pi)=\lambda \neq 0$ for $k \in \mathbb{Z}$.

(i) When $g'(k\pi)<0$, by theorem (\ref{thm:2.3}), there don't exist vortex solutions with finite energy which satisfy $\lim_{r \to \infty} h(r)=k\pi$.

(ii) When $g'(k\pi)>0$, for $G(x,k)=\frac{\lambda}{2}\sin^2 x \geq 0$, by theorem (\ref{thm:2.2}), there don't exist vortex solutions which satisfy $\lim_{r \to \infty} h(r)=k\pi$.

From above argument, we know there don't admit any vortex solutions with finite energy in this case.

\medskip
\noindent {\large \bf {Case four}}: $\lambda > 0,\,\omega >0$.

Now we have $g'(k\pi)=\lambda + \omega (-1)^{|k|}$ for $k \in \mathbb{Z}$.

(i) when $k=2l,l \in \mathbb{Z}$, it follows that
$$g'(2l\pi)=\lambda +\omega>0\ \text{and}\ G(x,2l)=\frac{\lambda}{2}\sin^2 x +\omega[1-\cos x]\geq 0.$$
By theorem (\ref{thm:2.2}), there don't exist vortex solutions which satisfy $\lim_{r \to \infty} h(r)=2l\pi$.

(ii) when $k=2l+1,l \in \mathbb{Z}$, it follows that $g'((2l+1)\pi)=\lambda -\omega$. For $\lambda<\omega$, by theorem (\ref{thm:2.3}), there don't exist vortex solutions with finite energy which satisfy $\lim_{r \to \infty} h(r)=(2l+1)\pi$.

From above argument, we deduce that if there exists a vortex solution with finite energy for $\lambda > 0,\omega >0$, there holds
$$0<\omega \leq \lambda\quad \text{and}\quad \lim_{r \to \infty}h(r)=(2l+1)\pi,\  \text{for some}\ l \in \mathbb{Z}.$$
Moreover, for $\omega<\lambda$, by lemma (\ref{lem:2.1}),  $h(r)$ converges to $(2l+1)\pi$ exponentially as $r \to \infty$.

\medskip
\noindent {\large \bf {Case five}}: $\lambda > 0,\,\omega <0$.

Following the similar argument as Case four, we have that if there exists a vortex solution with finite energy for $\lambda > 0,\omega <0$, there holds
$$-\lambda \leq \omega <0\quad \text{and}\quad \lim_{r \to \infty}h(r)=2l\pi,\  \text{for some}\ l \in \mathbb{Z}.$$
Moreover, for $-\lambda<\omega$, by lemma (\ref{lem:2.1}),  $h(r)$ converges to $2l\pi$ exponentially as $r \to \infty$.

\medskip
\noindent {\large \bf {Case six}}: $\lambda < 0,\,\omega >0$.

Now we have $g'(k\pi)=\lambda + \omega (-1)^{|k|}$ for $k \in \mathbb{Z}$.

(i)when $k=2l+1,l \in \mathbb{Z}$, it follows that $g'((2l+1)\pi)=\lambda -\omega<0$. By theorem (\ref{thm:2.3}), there don't exist vortex solutions with finite energy which satisfy $\lim_{r \to \infty} h(r)=(2l+1)\pi$.

(ii) when $k=2l,l \in \mathbb{Z}$, it follows that $g'(2l\pi)=\lambda +\omega$. For $0<\omega<-\lambda$, by theorem (\ref{thm:2.3}), there don't exist vortex solutions with finite energy which satisfy $\lim_{r \to \infty} h(r)=2l\pi$. For $\omega>-\lambda$, it's easy to verify that $G(x,2l)=\frac{\lambda}{2}\sin^2 x +\omega[1-\cos x]\geq 0$. Thus, by theorem (\ref{thm:2.2}), there don't exist vortex solutions which satisfy $\lim_{r \to \infty} h(r)=2l\pi$.

Now from above argument, if there exists a vortex solution with finite energy for $\lambda < 0,\omega >0$, there holds
$$-\lambda =\omega >0\quad \text{and}\quad \lim_{r \to \infty}h(r)=2l\pi,\  \text{for some}\ l \in \mathbb{Z}.$$

\medskip
\noindent {\large \bf {Case seven}}: $\lambda < 0,\,\omega <0$.

Following the similar argument as Case six, we have that if there exists a vortex solution with finite energy for $\lambda < 0,\omega <0$, there holds
$$\lambda =\omega <0\quad \text{and}\quad \lim_{r \to \infty}h(r)=(2l+1)\pi,\  \text{for some}\ l \in \mathbb{Z}.$$

Now, we complete the proof of theorem (\ref{thm:1}).

\medskip

{\bf Acknowledge} The author would like to thank his supervisors
Professor Weiyue Ding and Professor Youde Wang for their
encouragement and inspiring advices.

\noindent Ruiqi Jiang\\
Academy of Mathematics and Systems Science\\
Chinese Academy of Sciences,\\
Beijing 100190, P.R. China.\\
Email: jiangruiqi@amss.ac.cn

\end{document}